\documentclass[12pt]{article}
 \usepackage{amsmath,amsthm,amssymb,amscd}

\title{Flops and derived categories}
\author{Tom Bridgeland}

\date{}

\newcommand{\pro}[1]{\paragraph{\sc{(4.#1)}}}
\newtheorem{thm}{Theorem}[section]
\newtheorem{cor}[thm]{Corollary}
\newtheorem{prop}[thm]{Proposition}
\newtheorem{lemma}[thm]{Lemma}
\newenvironment{pf}{\paragraph{Proof}}{\qed\par\medskip}

\theoremstyle{definition}
\newtheorem{defn}[thm]{Definition}
\newtheorem{example}[thm]{Example}

\newcommand{\fs}{f_{S\,*}}

\newcommand{\Cp}{C'}
\newcommand{\HH}{\mathcal H}
\newcommand{\Phib}{{\Tilde{\Phi}}}

\newcommand{\isom}{\cong}

\newcommand{\tensor}{\otimes}
\newcommand{\Li}{\mathcal{L}}
\newcommand{\cl}{\operatorname{c}}
\newcommand{\st}{{s\,*}}
\newcommand{\PP}{\mathbb P}
\newcommand{\M}{\operatorname{\mathcal M}}
\newcommand{\C}{\mathbb C}
\newcommand{\pHilb}{\operatorname{P-Hilb}}
\newcommand{\Hilb}{\operatorname{Hilb}}
\newcommand{\PC}{\operatorname{{Per}}\,}
\newcommand{\Quot}{\operatorname{Quot}}

\newcommand{\I}{\mathcal I}
\newcommand{\E}{\mathcal E}
\newcommand{\F}{\mathcal F}
\newcommand{\G}{\mathcal G}
\newcommand{\W}{L}
\newcommand{\U}{\mathcal U}
\newcommand{\ts}{$t$-structure }
\newcommand{\tss}{$t$-structures }

\newcommand{\K}{K}
\newcommand{\A}{\mathcal A}
\newcommand{\B}{\mathcal B}
\newcommand{\CC}{\mathcal C}
\newcommand{\OO}{\mathcal O}
\newcommand{\into}{\hookrightarrow}
\newcommand{\onto}{\twoheadrightarrow}
\renewcommand{\P}{\mathcal P}

\newcommand{\id}{\operatorname{id}}

\newcommand{\D}{\operatorname{D}}
\newcommand{\Ext}{\operatorname{Ext}}
\newcommand{\Hom}{\operatorname{Hom}}

\newcommand{\SL}{\operatorname{SL}}
\newcommand{\eu}{\operatorname{\chi}}

\renewcommand{\L}{\mathbf L}
\newcommand{\R}{\mathbf R}
\newcommand{\Ltensor}{\stackrel{\mathbf L}{\tensor}}

\newcommand{\lRa}[1]{\xrightarrow{\ #1\ }}

\newcommand{\lra}{\longrightarrow}

\newcommand{\Bp}{\B\smash{^\perp}}

\makeatletter
 \newlength{\typesize}
\makeatother

\newlength{\vvoff}
\newlength{\hhoff}

\newcommand{\locateoffcenter}[1]{%
\addtolength{\vvoff}{-0.25\typesize}%
\raisebox{\vvoff}{\hspace{\hhoff}\makebox(0,0){\smash{#1}}}
}
\newcommand{\object}[1]{%
\setlength{\vvoff}{0pt}%
\setlength{\hhoff}{0pt}%
\locateoffcenter{#1}
}
\newcommand{\elabel}[1]{%
\setlength{\vvoff}{0.75\typesize}%
\setlength{\hhoff}{0pt}%
\locateoffcenter{#1}
}

\newcommand{\slabel}[1]{%
\setlength{\vvoff}{0pt}%
\setlength{\hhoff}{0.75\typesize}%
\locateoffcenter{#1}
}

\newcommand{\swlabel}[1]{%
\setlength{\vvoff}{-0.5\typesize}%
\setlength{\hhoff}{0.75\typesize}%
\locateoffcenter{#1}
}
\newcommand{\nwlabel}[1]{%
\setlength{\vvoff}{-0.5\typesize}%
\setlength{\hhoff}{-0.75\typesize}%
\locateoffcenter{#1}
}

\begin{document}
\maketitle



\section{Introduction}

This paper contains some applications
of Fourier-Mukai techniques to problems in birational geometry.
The main new idea is that flops occur naturally as
moduli spaces of perverse coherent sheaves.
As an application we prove

\begin{thm}
\label{main}
If $X$ is a projective threefold with terminal
singularities and
\begin{equation*}
 \setlength{\unitlength}{36pt}
 \begin{picture}(2.2,1.2)(0,0)
 \put(1,0){\object{$X$}}
 \put(0,1){\object{$Y_1$}}
 \put(2,1){\object{$Y_2$}}
 \put(0.25,0.75){\vector(1,-1){0.5}}
 \put(0.5,0.5){\nwlabel{$f_1$}}
 \put(1.75,0.75){\vector(-1,-1){0.5}}
 \put(1.5,0.5){\swlabel{$f_2$}}
 \end{picture}
\end{equation*}
are crepant resolutions,
then there is an equivalence of derived categories of coherent sheaves
$\D(Y_1)\lra \D(Y_2)$.
\end{thm}

The theorem implies in particular that birational Calabi-Yau threefolds
have equivalent derived categories and thus gives a new
proof of the theorem (due to V.V. Batyrev \cite{Ba})
that birational Calabi-Yau
threefolds have the same Hodge numbers. A.I. Bondal and D.O. Orlov
proved some special cases of Theorem \ref{main}
in \cite{BO}
and conjectured that the result held in general. Here we shall prove
Theorem \ref{main} using the by-now standard techniques of
Fourier-Mukai transforms, in particular the ideas developed in
\cite{BKR,Br2}.

\medskip

For simplicity, let
us suppose that $Y$ is a non-singular, projective
threefold and $f\colon Y\to X$ is a proper, birational morphism
contracting a single rational curve $C\isom\PP^1$ with
normal bundle 
\[\mathcal{N}_{C/Y}\isom \OO_C(-1)\oplus\OO_C(-1).\]

Using the theory of $t$-structures, we define an abelian category
$\PC(Y/X)\subset \D(Y)$ whose objects we call \emph{perverse}
(or \emph{perverse coherent})
sheaves on $Y$. A short exact sequence in $\PC(Y/X)$
is just a triangle in $\D(Y)$ whose vertices are all
objects of $\PC(Y/X)$.
The next step is to construct moduli spaces of perverse sheaves.
To do this we introduce a stability condition.
A \emph{perverse point sheaf} is then defined to be
a stable perverse sheaf which has the same numerical
invariants as the structure sheaf of a
point of $Y$.

Structure sheaves of points $y\in Y$
are objects of the
category $\PC(Y/X)$, and are stable for $y\in Y\setminus C$.
For $y\in C$, the sheaf
$\OO_y$ fits into the exact sequence
\begin{equation}
\label{a}
0\lra \OO_C(-1)\lra \OO_C\lra \OO_y\lra 0.
\end{equation}
It turns out that $\OO_C$ is a perverse sheaf, but $\OO_C(-1)$ is
not, so that the triangle in $\D(Y)$ arising from (\ref{a})
does not define an exact sequence in
$\PC(Y/X)$. However the complex obtained by shifting $\OO_C(-1)$
to the left by one place \emph{is} a perverse sheaf, so there
is an exact sequence of perverse sheaves
\begin{equation}
\label{b}
0\lra \OO_C\lra\OO_y\lra\OO_C(-1)[1]\lra 0,
\end{equation}
which should be thought of as destabilizing $\OO_y$.

Flipping the extension of perverse sheaves (\ref{b}) gives
stable objects of $\PC(Y/X)$ fitting into an exact sequence of
perverse sheaves
\begin{equation}
\label{c}
0\lra \OO_C(-1)[1]\lra E\lra \OO_C\lra 0.
\end{equation}
These perverse point sheaves $E$ are
not sheaves, indeed any such object has two non-zero
homology sheaves $H_1(E)=\OO_C(-1)$ and $H_0(E)=\OO_C$.
We shall use geometric invariant theory
to construct a fine moduli space $W$ parameterizing perverse point
sheaves on $X$. Roughly speaking, the space $W$ is obtained from $X$
by replacing the rational curve $C$ parameterising extensions
(\ref{b}) by another rational curve $C'$ parameterising
extensions (\ref{c}).

The push-down $\R f_*(E)$ of a
perverse point sheaf $E$ is always the structure sheaf of a point $x\in X$,
so there is a natural map $g\colon W\to X$.
Moreover, the general point of $W$ corresponds to
the structure sheaf of a point $y\in Y\setminus C$, so $g$ is
birational. Thus there is a diagram of birational morphisms
\begin{equation*}
 \setlength{\unitlength}{36pt}
 \begin{picture}(2.2,1.2)(0,0)
 \put(1,0){\object{$X$}}
 \put(0,1){\object{$W$}}
 \put(2,1){\object{$Y$}}
 \put(0.25,0.75){\vector(1,-1){0.5}}
 \put(0.5,0.5){\nwlabel{$g$}}
 \put(1.75,0.75){\vector(-1,-1){0.5}}
 \put(1.5,0.5){\swlabel{$f$}}
 \end{picture}
\end{equation*}

The techniques developed in \cite{BKR,Br2}
allow us to use the intersection theorem to
show that $W$ is non-singular, and that
the universal family of perverse sheaves on
$W\times Y$ induces a Fourier-Mukai transform $\D(W)\lra\D(Y)$.
An easy argument then shows that
$g\colon W\to X$ is the flop of $f\colon Y\to X$.
Theorem \ref{main} follows from this because crepant resolutions
of a terminal threefold are related by a finite chain of flops (see
 \cite{Ko}).

\medskip

The hard work in this paper goes into constructing the moduli space
of perverse point sheaves. It turns out that
the correct stability condition to impose on these objects
is that they should be quotients of $\OO_Y$ in the category $\PC(Y/X)$.
Thus each perverse point sheaf $E$ fits into an exact sequence
of perverse sheaves
\begin{equation}
0\lra F\lra \OO_Y\lra E\lra 0
\end{equation}
and the space $W$ is really a sort of perverse Hilbert scheme
parameterising
perverse quotients of $\OO_Y$. The corresponding
subobjects $F\subset\OO_Y$ are simple, rank one sheaves,
in general with torsion. Thus in the first place we construct a 
moduli space of simple sheaves $F$ on $Y$ and then use this space to
parameterise the corresponding perverse point sheaves $E$.

The theory of perverse
sheaves developed below is valid for any small contraction of
canonical threefolds $f\colon Y\to X$.
It seems natural to speculate that when $-\K_Y$ is $f$-ample
the resulting moduli space of perverse sheaves $W$ is the \emph{flip} of
$f$. In that case one would not expect
a derived equivalence between $Y$ and $W$, but rather an
embedding of derived categories $\D(W)\into\D(Y)$.
What prevents us from proving
such a result is our inability to do Fourier-Mukai on singular spaces.
There is some hope that a better understanding of the
mathematics surrounding the intersection theorem might allow flips
to be studied in this way. This would be interesting for several reasons,
not least because it would give a simpler and more conceptual
proof of the existence of threefold flips. For now, however, this
remains pure speculation!

\medskip

The plan of the paper is as follows. Section 2 contains the
basic definitions we need from the theory of triangulated categories.
In Section 3 we define the category of perverse coherent sheaves
and derive some of its basic properties. We also state
Theorem \ref{rep} which guarantees the existence of
fine moduli spaces of
perverse point sheaves. In Section 4 we assume this result and use
it to prove Theorem
\ref{main}.  The proof of Theorem \ref{rep} is given in Sections 5 and
6.

\medskip

\noindent {\bf Notation.} All schemes $X$ are assumed to be of finite type over
$\C$ and all points are closed points. $\D(X)$ denotes the unbounded derived
category of coherent sheaves throughout. More precisely $\D(X)$ is the
subcategory of the derived category of quasi-coherent $\OO_X$-modules
consisting of complexes with coherent cohomology sheaves. The
full subcategory of complexes with bounded cohomology sheaves is
denoted $\D^b(X)$. The $i$th cohomology sheaf of an
object $E\in\D(X)$ is denoted $H^i(E)$ and the $i$th homology sheaf
by $H_i(E)$. Thus
 $H_i(E)=H^{-i}(E)$.


\section{Admissible subcategories}

This section contains some ideas from the general theory of triangulated
categories. In particular, we define
semi-orthogonal decompositions and $t$-structures.
In fact, the first of
these concepts is a special case of the second, but we give the definitions
separately, since one tends to
think of the two structures rather differently. $t$-structures were
introduced in \cite{BBD} in order to
define perverse sheaves on stratified spaces. Semi-orthogonal
decompositions also appear in \cite{BBD} but their geometrical
significance was first properly exploited by Bondal and Orlov
\cite{Bo,BO}.
We fix a triangulated
category $\A$ throughout, with its shift functor
$T\colon\A\to\A\colon a\mapsto a[1]$.

\medskip

In the context of birational geometry, the key point to note about
derived categories is that performing a contraction
corresponds to passing to a triangulated subcategory.
More specifically one should consider so-called
\emph{admissible subcategories}.

\begin{defn}
A \emph{right admissible subcategory} of $\A$ is a full subcategory
$\B\subset\A$ such that the inclusion
functor $\B\into\A$ has a right adjoint.
\end{defn}

\noindent Given a full subcategory $\B\subset\A$
one defines the right orthogonal $\Bp\subset\A$
to be the full subcategory
\[\Bp=\{a\in\A : \Hom_{\A}(b,a)=0 \text{ for all } b\in \B\}.\]
One can easily show \cite{Bo} that
if a full subcategory $\B\subset\A$ is right admissible
then every object $a\in
A$ fits into a triangle
\[b\lra a\lra c\lra b[1]\]
with $b\in\B$ and $c\in\Bp$.

\begin{defn}
A \emph{triangulated subcategory} of $\A$ is a
full subcategory $\B\subset\A$
which is closed under shifts, that is $\B[1]=\B$, such that
any triangle
$b_1\lra b_2\lra c\lra b_1[1]$
in $\A$ with $b_1,b_2\in\B$ has $c\in B$ also.
\end{defn}

\noindent Clearly the right orthogonal of a triangulated category is itself
triangulated.
If a triangulated subcategory $\B\subset\A$ is right admissible
we say that $\A$ has a \emph{semi-orthogonal decomposition}
into the subcategories $(\Bp,\B)$; one should think of $\A$ as being
built up from these two smaller triangulated categories. Important
examples of semi-orthogonal decompositions are given by the
following result.

\begin{prop}
\label{first}
Let $f\colon Y\to X$ be a morphism of projective varieties
such that $\R f_* (\OO_Y)=\OO_X$. Then the functor
\[ \L f^*\colon \D(X)\lra \D(Y) \]
embeds $\D(X)$ as a right admissible
triangulated subcategory of $\D(Y)$.
\end{prop}

\begin{pf}
The functor $\L f^*$ has the right adjoint $\R f_*$ and the composite
$\R f_* \circ \L f^*$ is the identity on $\D(X)$ by the projection
formula and the assumption that $\R f_* (\OO_Y)=\OO_X$. It follows
that $\L f^*$ is fully faithful so $\D(X)$ can be identified with its
image under $\L f^*$. Note that if $X$ is non-singular then
$\L f^*$ embeds $\D^b(X)$ in $\D^b(Y)$, but that this is no longer
true when we allow singularities.
\end{pf}

The Grauert-Riemenschneider vanishing theorem shows that the
hypotheses of Proposition \ref{first} hold whenever $f\colon Y\to X$
is a morphism of projective varieties such
that $Y$ has rational singularities and $-\K_Y$ is $f$-ample.

\begin{cor}
Let $f\colon Y\to X$ be an extremal contraction of a canonical
threefold. Then $\D(X)$ is a right admissible triangulated
subcategory of $\D(Y)$.
\qed
\end{cor}

As we mentioned in the introduction, it is possible
that flips also induce embeddings of derived categories.
If this were true, one would be able to interpret the action of
the minimal model program on a variety $X$ as picking out
some minimal admissible subcategory of $\D(X)$. 

\medskip

Recall that an abelian category $A$ sits inside its derived category
$\D(A)$ as the subcategory of complexes
whose cohomology is concentrated in degree zero.
There are by now plenty of examples of interesting algebraic and
geometrical relationships which can be described by an equivalence of
derived categories $\D(A)\lra\D(B)$. Such equivalences will usually
not arise from an equivalence of the underlying abelian categories $A$
and $B$, indeed, this is why one must use derived categories.
Changing perspective slightly one could think of a derived equivalence
as being described by a single triangulated
category with two different abelian categories sitting inside it.
The theory of \tss is the tool which
allows one to see these different categories.

\begin{defn}
A \emph{\ts}on $\A$ is a right admissible
subcategory $\A^{\leq 0}\subset\A$ which is preserved by left shifts,
that is $\A^{\leq 0}\,[1]\subset \A^{\leq 0}$.
\end{defn}

\noindent Given a \ts $\A^{\leq 0}$ on $\A$ one defines
$\A^{\leq i}=\A^{\leq 0}[-i]$ and
$\A^{\geq i}=(\A^{\leq {i-1}})\smash{^\perp}$. One also writes
$\A^{<i}=\A^{\leq i-1}$ and $\A^{>i}=\A^{\geq i+1}$.

\begin{defn}
The \emph{heart}
(or \emph{core}) of the \ts $\A^{\leq 0}\subset\A$
is the full subcategory
$\HH=\A^{\leq 0}\cap\A^{\geq 0}$.
\end{defn}

\noindent It was proved in \cite{BBD}
that the heart of a $t$-structure is an abelian category.
Short exact sequences
$0\lra a_1\lra a_2\lra a_3\lra 0$
in $\HH$ are determined by triangles
$a_1\lra a_2\lra a_3\lra a_1[1]$
in $\A$ with $a_i\in\HH$ for all $i$.

\smallskip

The basic example is the \emph{standard $t$-structure} on
the derived category $\D(A)$ of an abelian category $A$, given by
\begin{eqnarray*}
\A^{\leq 0}&=&\{E\in \D(A): H^i(E)=0 \text{ for all } i>0\}, \\
\A^{\geq 0}&=&\{E\in \D(A): H^i(E)=0 \text{ for all } i<0\}.
\end{eqnarray*}
The heart is the original abelian category $A$.
To give another example, suppose
that $\D(A)\lra\D(B)$ is an equivalence of derived categories.
Then pulling back the standard \ts on $\D(B)$
gives a \ts on $\D(A)$ whose heart is the abelian
category $B$.

Further examples are provided by admissible triangulated
subcategories $\B\subset\A$. Any such subcategory defines a \ts
whose heart is trivial. In fact the converse is true:
a $t$-structure $\A^{\leq
0}\subset\A$ satisfying  $\A^{\leq 0}\cap\A^{\geq 0}=0$
is actually a triangulated subcategory of $\A$.
We shall not need this fact and the proof is left to the reader.


\section{Perverse coherent sheaves}

In this section we define the category of perverse sheaves with which
we shall be working for the rest of the paper.
It is objects of this category which will be naturally parameterised by
the points of a flop.
Let  $f\colon Y\to X$ be a birational morphism of
projective varieties. We shall make two assumptions, firstly
that $\R f_* \OO_Y=\OO_X$, and secondly that $f$ has relative
dimension one.  The example we have in mind is a small
contraction of a canonical threefold.

\medskip

Let us write $\A=\D(Y)$ and $\B=\D(X)$.
By Proposition \ref{first}, we may identify
$\B$ with a right admissible
triangulated subcategory of $\A$. Thus there is a semi-orthogonal
decomposition $(\CC,\B)$ where
\[\CC=\Bp=\{E\in\D(Y) : \R f_* (E) =0 \}.\]
Note that objects of $\CC$ are supported on the exceptional locus of
$f$.

\begin{lemma}
\label{easy}
An object $E\in\D(Y)$ lies in $\CC$ precisely when  its cohomology
sheaves $H^i(E)$ lie in $\CC$.
\end{lemma}

\begin{pf}
There is a spectral sequence $\R^p f_* H^q(E) \implies H^{p+q} \,\R f_* (E)$
which degenerates because $f$ has relative dimension one.
\end{pf}

Note that the functor $\R f_*$ has the left adjoint
$\L f^*$ and the right adjoint $f^!$. In this situation
one may obtain $t$-structures on $\A$ by glueing
$t$-structures on $\B$ and $\CC$. For details see 
\cite[1.4.8 - 10]{BBD} or \cite[Ex. IV.4.2 (c)]{GM}.
Lemma \ref{easy} allows one to use the standard \ts on $\A$ to induce a
\ts $\CC^{\leq 0}=\CC\cap\A^{\leq 0}$ on $\CC$ in the obvious way.
Shifting this by an integer $p$ and glueing it
to the standard $t$-structure on $\B=\D(X)$
gives a \ts on $\A$ satisfying
\begin{eqnarray*}
^p\!\A^{\leq 0}&=&\{E\in \A:\R f_* (E)\in \B^{\leq 0} \text{ and }
\Hom_{\A}(E,C)=0\text{ for all }C\in\CC^{>p} \}, \\
^p\!\A^{\geq 0}&=&\{E\in \A:\R f_* (E)\in \B^{\geq 0} \text{ and }
\Hom_{\A}(C,E)=0\text{ for all }C\in\CC^{<p} \}.
\end{eqnarray*}
The heart of this $t$-structure is the abelian category
\[\smash{^p}\PC(Y/X)=\,^p\!\A^{\leq 0}\cap\,^p\!\A^{\geq 0}.\]
The integer $p$ should be thought of as a choice of perversity.
We shall be mainly interested in the case $p=\smash{-}1$, and we refer
to objects of the category
\[\PC(Y/X)=\smash{^{\smash{-}1}}\PC(Y/X)\]
as \emph{perverse} (or \emph{perverse coherent}) sheaves.
The lemma below gives an explicit description of this category.

\begin{lemma}
\label{perverse}
An object $E$ of $\D(Y)$ is a perverse sheaf if and only if the
following three conditions are satisfied:
\begin{description}
\item[(a)] $H_i(E)=0$ unless $i=0$ or $1$,

\item[(b)] $\R^1 f_*\, H_0(E)=0$ and $\R ^0 f_*\, H_1(E)=0$,

\item[(c)] $\Hom_X(H_0(E),C)=0$ for any sheaf $C$ on $Y$ satisfying $\R f_*(C)=0$.
\end{description}
\end{lemma}

\begin{pf}
Suppose $E$ is a perverse sheaf.
The condition that  $\R f_*(E)$ is a sheaf on $X$,
together with the spectral sequence of Lemma \ref{easy}, gives
condition (b)
and implies that $\R f_*\, H_i(E)=0$ unless $i=0$ or $1$.

Let $\tau_{< i}$ and $\tau_{> i}$ be the truncation functors
of the standard $t$-structure on $\D(Y)$.
There are natural maps $\tau_{< i}\,E\to E$ and $E\to\tau_{> i}\,E$.
Then $\tau_{> 0}\,E=0$
because $\tau_{>0}\,E\in\CC^{>0}$.
Similarly $\tau_{<-1}\,E=0$ because $\tau_{<-1}\,E\in\CC^{<-1}$.
This proves condition (a). Condition (c) is clear, since any
non-zero map from $H_0(E)$ to a sheaf $C\in\CC^{>-1}$ induces a non-zero
morphism $E\to C$ in $\D(Y)$.

The converse is easy and is left to the reader.
\end{pf}

\begin{defn}
We shall say that two objects $A_1$ and $A_2$ of $\D^b(Y)$
are \emph{numerically equivalent} if for any locally-free sheaf
$L$ on $Y$ one has $\eu(L,A_1)=\eu(L,A_2)$.
\end{defn}

\noindent Recall that for objects $L$ and $A$ of $\D^b(Y)$
with $L$ of finite homological dimension
\[\eu(L,A)=\sum_i (-1)^i \dim \Hom^i_{\D(Y)}(L,A).\]
Thus if $Y$ is a non-singular projective variety, then by the
Riemann-Roch theorem, two objects of $\D^b(Y)$ are
numerically equivalent precisely when they have the same Chern character.

\begin{defn}
An object $F$ of $\D(Y)$ is a \emph{perverse ideal sheaf}
if there is an injection $F\into\OO_Y$ in the category $\PC(Y/X)$. An object
$E$ of $\D(Y)$ is a \emph{perverse structure sheaf}
if there is a surjection $\OO_Y\onto E$ in the category $\PC(Y/X)$. A
\emph{perverse point sheaf} is a perverse structure sheaf which is numerically equivalent to the
structure sheaf of a point $y\in Y$.
\end{defn}

Thus a perverse ideal sheaf $F$ determines and is determined by a perverse
structure sheaf $E$, which fit together in an exact sequence
of perverse sheaves
\begin{equation}
\label{birthday}
0\lra F\lra \OO_Y\lra E\lra 0.
\end{equation}
Applying the cohomology functor to the above exact sequence
shows that perverse ideal sheaves are
actually sheaves, that is satisfy $H_i(F)=0$ for $i\neq 0$.

\begin{example}
Let us suppose, as in the introduction,
that $f\colon Y\to X$ is the contraction of a
non-singular rational curve $C$ with normal bundle
$\OO_C(-1)\oplus\OO_C(-1)$ on a non-singular projective
threefold $Y$. For any point $y\in C$ there is an exact sequence
of sheaves
\[0\lra\I_C\lra \I_y\lra \OO_C(-1)\lra 0.\]
This shows that $\I_y$ is not a perverse sheaf, and it follows that
although $\OO_y$ is a perverse sheaf, it is not a quotient of $\OO_Y$
in $\PC(Y/X)$. Consider instead non-trivial extensions of the form
\[0\lra\OO_C(-1)\lra F\lra\I_C\lra 0.\]
One can easily calculate that $\Ext^1_Y(\I_C,\OO_C(-1))=\C^2$, so the
set of such sheaves $F$ is parameterised by a rational curve.
Composing the map $F\to\I_C$ with the inclusion $\I_C\subset\OO_Y$
gives a non-zero morphism $F\to\OO_Y$ and we take $E$ to be its cone.
In this way we obtain an exact sequence of perverse sheaves
\[0\lra F\lra \OO_Y\lra E\lra 0\]
with $H_1(E)=\OO_C(-1)$ and $H_0(E)=\OO_C$. Thus $E$ is a perverse
point sheaf.

The flop of $Y$ along $C$ is a non-singular threefold $W$ with a
morphism $g\colon W\to X$ contracting a single rational curve $\Cp$. We
shall show that the points of $W$ parameterise perverse point
sheaves on $Y$. The perverse point corresponding to
a point $w\in W\setminus \Cp$ 
is a point $y\in Y\setminus C$, whereas
the points of $\Cp$ correspond to the perverse point sheaves
$E$ described above.
\end{example}

\smallskip

\begin{lemma}
\label{nearly}
Let $E_1$ and $E_2$ be perverse point sheaves on $Y$. Then
 \[
 \Hom_{\D(Y)}(E_1,E_2)
 =\begin{cases}
 \,\C &\text{if $E_1=E_2$,} \\
 \, 0 &\text{otherwise.}
\end{cases}
 \]
\end{lemma}

\begin{pf}
If $E$ is a perverse point sheaf then $\R f_*(E)$ is the structure
sheaf of a point of $X$, so $\Hom_{\D(Y)}(\OO_Y,E)=\C$. Taking
$\Hom$s of the exact sequence
(\ref{birthday}) into $E$ shows that $\Hom_{\D(Y)}(E,E)=\C$.

Suppose there is a non-zero morphism $\theta\colon E_1\to E_2$.
Taking $\Hom$s of (\ref{birthday}) into $E_2$ shows that the unique
map $\OO_Y\to E_2$ must factor via $\theta$. In particular, $\theta$ is
surjective in $\PC(Y/X)$. But then the kernel $K$ of $\theta$ in $\PC(Y/X)$
is numerically equivalent to zero, and this implies $K=0$,
so $\theta$ is an isomorphism.
\end{pf}

\medskip

Let $S$ be a scheme. Given a point $s\in S$, let
$j_s\colon \{s\}\times Y\into S\times Y$
 be the embedding. A family of sheaves
on $Y$ over $S$ is just
an object $\F$ of $\D(S\times Y)$ such that for each point $s\in S$ the
object $\F_s=\L j_s^* (\F)$ of $\D(Y)$ is a sheaf.
Indeed, by \cite[Lemma 4.3]{Br1}, this condition implies that
$\F$ is actually a
sheaf on $S\times Y$, flat over $S$, so that $\F$ defines
a family of sheaves in the usual sense.
Once this observation has been made 
it is clear what the correct definition of a family of perverse
sheaves should be.

\begin{defn}
A family of perverse sheaves
on $Y$ over a scheme $S$ is
an object $\E$ of $\D(S\times Y)$ such that for each point $s\in S$ the
object $\E_s=\L j_s^* (\E)$ of $\D(Y)$ is a perverse sheaf. Two such
families $\E_1$ and $\E_2$ are equivalent if $\E_2=\E_1\tensor L$ for
some line bundle $L$ pulled back from $S$.
\end{defn}

The proof of the following theorem will be given in Sections 5 and 6 below.

\begin{thm}
\label{rep}
The functor which assigns to a scheme $S$ the set of equivalence
classes of families of
perverse point sheaves on $Y$ over $S$ is representable by a projective
scheme $\M(Y/X)$.
\end{thm}

We conclude this section with the
following base-change result.

\begin{prop}
\label{basechange}
Let $S$ be a scheme and $\E$ a family of perverse sheaves on $Y$
over $S$. Put $f_S=\id_S\times f$.
Then $\G=\fs (\E)$ is a family of sheaves
on $X$ over $S$, and for any point $s\in S$ there is an isomorphism
of sheaves $\G_s=\R f_*(\E_s)$.
\end{prop}

\begin{pf}
Fix a point $s\in S$ and consider the diagram
\begin{equation*}
\setlength{\unitlength}{48pt}
\begin{picture}(2.2,1.2)(0,0)
\put(0,0){\object{$X$}}
\put(2,0){\object{$S\times X$}}
\put(0,1){\object{$Y$}}
\put(2,1){\object{$S\times Y$}}
\put(0.3,0){\vector(1,0){1.2}}
\put(0,0.8){\vector(0,-1){.6}}
\put(0.3,1){\vector(1,0){1.2}}
\put(2,0.8){\vector(0,-1){.6}}
\put(1,0){\elabel{$i_s$}}
\put(0,.55){\slabel{$f$}}
\put(1,1){\elabel{$j_s$}}
\put(2,.55){\slabel{$f_S$}}
\end{picture}
\end{equation*}
where $i_s\colon \{s\}\times X\into S\times X$
is the embedding.

If $p\colon S\times X\to S$ is the projection map, flat
base-change shows that $i_\st \OO_X=p^* \OO_s$.
It follows that $\L f_S^* \circ i_\st \OO_X = j_\st \OO_Y$.
The projection formula gives isomorphisms
\begin{eqnarray*}
 i_\st\circ \R f_*(\E_s) &=& \R \fs\circ j_\st
\circ \L j_s^* (\E) =  \R \fs (j_\st \OO_Y \Ltensor\E) \\
&=& i_\st (\OO_X) \Ltensor \R \fs (\E)= i_\st\circ \L i_s^*(\G).
\end{eqnarray*}
The functor $i_\st$ is exact and fully faithful on the category of
sheaves on $X$.
Thus $\R f_*(\E_s)$ is a sheaf precisely when $\L i_s^*(\G)$ is.
Since $\E$ is a family of perverse sheaves,
 $\R f_*(\E_s)$ is a sheaf for all $s\in S$,
so $\G$ is a sheaf, flat over $S$. The result follows.
\end{pf}


\section{Flops and the derived equivalence}

In this section we shall show how fine moduli spaces of perverse point
sheaves give rise to Fourier-Mukai type equivalences of derived
categories. To do this we shall assume that fine moduli spaces
of perverse point sheaves exist (as in Theorem \ref{rep}) and apply
the techniques of \cite{BKR,Br2}. In this way we obtain a proof
of Theorem \ref{main}.

\pro{1} Let $X$ be a projective threefold with
Gorenstein, terminal singularities. Recall that a \emph{crepant
resolution}
is a morphism $f\colon Y\to X$ from a non-singular projective variety
$Y$, such that $f^*\omega_X=\omega_Y$. Any such resolution
satisfies $\R f_*(\OO_Y)=\OO_X$ and
contracts only a finite number of curves. Thus
the open subset $U\subset X$ over which $f$ is an
isomorphism is the complement of a finite set of points.

By Theorem \ref{rep} there is a fine moduli space $\M(Y/X)$ of perverse
point sheaves on $Y$.
Each point $y\in f^{-1}(U)$ is a perverse point sheaf so there
is an embedding $U\into \M(Y/X)$.
Let $W\subset\M(Y/X)$ be the
irreducible component of $\M(Y/X)$ containing the image of this
morphism. In fact, it is possible to prove,
as in \cite[Section 8]{BKR},
that $\M(Y/X)$ is
irreducible, so that $W=\M(Y/X)$, but we shall not need this.

Let $\P$ be a universal object on $W\times Y$.
Thus $\P$ is an object of $\D(W\times Y)$ such that
the perverse point sheaf on $Y$
corresponding to a point $w\in W$ is the object $\P_w=\L i_w^*(\P)$,
where $i_w\colon \{w\}\times Y\into W\times Y$ is the embedding.

\pro{2} By Proposition \ref{basechange}, the sheaf $\R(\id_W\times f)_*(\P)$
is a family of structure sheaves of points on $X$ over $W$, and
therefore, up to a twist by the pullback of a line bundle from $W$, is
the structure sheaf of the graph $\Gamma(g)\subset W\times X$
of some morphism $g\colon W\to X$. Thus twisting $\P$ by the pullback
of a line bundle from $W$, we can assume that
\begin{equation}
\label{iso}
\R(\id_W\times f)_*(\P)=\OO_{\Gamma(g)}.
\end{equation}
With this condition $\P$ is uniquely defined.
The morphism $g$ is birational because for any point $x\in U$ there is
only one object $E$ of $\D(Y)$ satisfying $\R f_*(E)=\OO_x$.
Thus there is a diagram of birational morphisms
\begin{equation*}
 \setlength{\unitlength}{36pt}
 \begin{picture}(2.2,1.4)(0,0)
 \put(1,0){\object{$X$}}
 \put(0,1){\object{$W$}}
 \put(2,1){\object{$Y$}}
 \put(0.25,0.75){\vector(1,-1){0.5}}
 \put(0.5,0.5){\nwlabel{$g$}}
 \put(1.75,0.75){\vector(-1,-1){0.5}}
 \put(1.5,0.5){\swlabel{$f$}}
 \end{picture}
\end{equation*}

\setcounter{thm}{2}

\pro{3} The scheme $W$ is a non-singular projective variety
and $g\colon W\to X$ is a crepant resolution. Furthermore,
the Fourier-Mukai functor
\[\Phi(-)=\R\pi_{Y,*}(\P\Ltensor\pi_{W}^*(-))\colon\D(W)\lra\D(Y),\]
is an equivalence of categories
which takes $\D^b(W)$ into $\D^b(Y)$.

\begin{pf}Each object $\P_w$ has bounded homology sheaves, and $Y$ is
non-singular, so the object $\P$ has finite homological dimension.
It follows that the functor $\Phi$ takes $\D^b(W)$ into $\D^b(Y)$.

For each point $w\in W$ the object $\P_w$ is simple,
so its support is connected, and
since $\R f_* (\P_w)=\OO_x$, where $x=g(w)$, it follows that
$\P_w$ is supported on the fibre of $f$ over $x$.
Since $f$ is crepant this implies that
$\P_w\tensor\omega_Y=\P_w$.

Given distinct points
$w_1, w_2\in W$, Serre duality together with Lemma \ref{nearly} shows that
\[
\Hom_{\D(Y)}^i(\P_{w_1},\P_{w_2})=0
\]
unless $g(w_1)=g(w_2)$ and $1\leq i\leq 2$.
The argument of \cite[Section 6]{BKR} then implies that $W$ is non-singular,
$g$ is crepant and $\Phi$ is an equivalence.
\end{pf}

\pro{4} An immediate consequence
of the isomorphism (\ref{iso}) is that there is a commutative diagram
of functors
\begin{equation*}
 \setlength{\unitlength}{48pt}
\begin{picture}(2.2,1.4)(0,0)
 \put(1,0){\object{$\D(X)$}}
 \put(0,1){\object{$\D(W)$}}
 \put(2,1){\object{$\D(Y)$}}
 \put(0.25,0.75){\vector(1,-1){0.5}}
 \put(0.4,0.5){\nwlabel{$\R g^*$}}
 \put(1.75,0.75){\vector(-1,-1){0.5}}
 \put(1.6,0.5){\swlabel{$\R f^*$}}
 \put(0.5,1){\vector(1,0){1}}
 \put(1,1){\elabel{$\Phi$}}
 \end{picture}
\end{equation*}

\pro{5}If $C$ is a sheaf on $W$ satisfying $\R g_*(C)=0$ then
$\Phi(C)[-1]$ is
a sheaf on $Y$.

\begin{pf}
First suppose that $A$ is an object of $\D(W)$ with $\R g_*(A)=0$
and such that $B=\Phi(A)[-1]$ is a sheaf on $Y$.
By (4.4) one has $\R f_*(B)=0$ so by Lemma \ref{perverse} the object
$B[1]$ is a perverse sheaf on $Y$. Thus for any point
$w\in W$,
\[\Hom_{\D(W)}^i(A,\OO_w)
=\Hom_{\D(Y)}^i(B[1],\P_w)=0\text{ unless }0\leq i\leq 3.\]
Moreover, since $\P_w$ is a quotient of $\OO_{Y}$ in $\PC(Y/X)$,
\[\Hom_{\D(Y)}^3(B[1],\P_w)=\Hom_{\D(Y)}^0(\P_w,B[1])=0.\]
Thus the object $A$ has homological dimension at most two, and is
supported in codimension at least two. It
follows from this that $A$ is a sheaf on $W$
(see, for example, \cite[Lemma 4.2]{Br2}).

Now assume that $C$ is a sheaf on $W$ satisfying $\R g_*(C)=0$
and suppose that $D=\Phi(C)[-1]$ is not a sheaf on $Y$.
As in the proof of
Lemma \ref{perverse}, we can find a sheaf $B$ on $Y$ satisfying $\R
f_*(B)=0$, and an integer $i<0$, such that
one of $\Hom_{\D(Y)}^i(B,D)$ or $\Hom_{\D(Y)}^i(D,B)$ is non-zero.
Since $\Phi$ is
an equivalence, $B=\Phi(A)[-1]$ for some object $A$ of $\D(W)$, and by
the first part $A$ is a sheaf. This implies
that one of the spaces
$\Hom_{\D(W)}^i(A,C)$ or $\Hom_{\D(W)}^i(C,A)$ is non-zero, which
is impossible since $A$ and $C$ are both sheaves.
\end{pf}

\pro{6}The variety $W=\M(Y/X)$ is the flop of $f\colon Y\to X$,
that is, if $D$ is a divisor on $W$ such that $-D$ is $g$-nef, then
its proper transform $D'$ on $Y$ is $f$-nef.

\begin{pf}Let $C$ be
a rational curve on $W$ contracted by $g$.
Put $M=\Phi(\OO_W(D))$ and $N=\Phi(\OO_C(-1))$. Then
by Riemann-Roch,
\[\eu(M,N)=\eu(\OO_W(D),\OO_C(-1))=-D\cdot C\geq 0.\]
Over the open subset $f^{-1}(U)$ of $Y$, $M$
is isomorphic to $\OO_Y(D')$, and it follows that $\cl_1(M)=[D']$.
By (4.5) the object $N[-1]$ is a sheaf
supported on some curve $C'$ of $Y$ which is contracted by $f$. It
follows that $\eu(M,N)=D'\cdot C'$
and hence the result.
\end{pf}

\pro{7}Define full triangulated subcategories
\[\CC(W/X)\into \D(X) \qquad \CC(Y/X)\into \D(Y)\]
consisting of objects
satisfying $\R g_*(C)=0$ and $\R f_*(C)=0$ respectively.
These categories inherit $t$-structures from the standard
$t$-structures on $\D(X)$ and $\D(Y)$, as in Section 3.
There is a commutative diagram of functors
\begin{equation*}
 \setlength{\unitlength}{48pt}
 \begin{picture}(4.2,1.4)(0,0)
 \put(0.15,0){\object{$\CC(Y/X)$}}
 \put(0.15,1){\object{$\CC(W/X)$}}
 \put(4,0){\object{$\D(X)$}}
 \put(4,1){\object{$\D(X)$}}
 \put(2,0){\object{$\D(Y)$}}
 \put(2,1){\object{$\D(W)$}}
 \put(0.7,0){\vector(1,0){.8}}
 \put(2.5,0){\vector(1,0){1}}
 \put(0.7,1){\vector(1,0){.8}}
 \put(2.5,1){\vector(1,0){1}}
 \put(3,1){\elabel{$\R g_*$}}
 \put(3,0){\elabel{$\R f_*$}}
 \put(0.15,.8){\vector(0,-1){0.6}}
 \put(0.2,.5){\slabel{$\Phib$}}
 \put(2,.8){\vector(0,-1){0.6}}
 \put(2.05,.5){\slabel{$\Phi$}}
 \put(4,.8){\vector(0,-1){0.6}}
 \put(4.05,.5){\slabel{$\id$}}
 \end{picture}
\end{equation*}
\smallskip
in which the rows are exact sequences of triangulated
categories. By (4.5), the equivalence $\Phib[-1]$
is $t$-exact, that is preserves the $t$-structures.

\pro{8}
There is a chain
of exact equivalences of abelian categories
\[\cdots\lra
  \smash{^{-1}}\PC(W/X)\lra
  \smash{^0}\PC(Y/X)\lra
  \smash{^1}\PC(W/X)\lra
  \smash{^2}\PC(Y/X)\lra\cdots.\]
Indeed, it follows from (4.7) that
for any integer $p$ the functor $\Phi$ induces an exact equivalence
\[\smash{^p}\PC(W/X)\isom\smash{^{p+1}}\PC(Y/X),\]
and since the flopping operation  is an involution we may interchange
$Y$ and $W$.


\section{Perverse ideal sheaves}

In this section we use
geometric invariant theory to construct fine moduli spaces
of perverse ideal sheaves.
As in Section 3, let $f\colon Y\to X$ be a birational morphism
of projective varieties of relative dimension one and satisfying
$\R f_*(\OO_Y)=\OO_X$. Our first task is to identify which objects of
$\D(Y)$ are perverse ideal sheaves.

\begin{prop}
\label{id}
A perverse ideal sheaf on $Y$ is, in particular, a sheaf on $Y$.
Furthermore, a sheaf on $Y$ is a perverse ideal sheaf if and only if
the following two conditions are satisfied:
\begin{description}
\item[(a)]  the sheaf $f_*(F)$ on $X$ is an ideal sheaf,

\item[(b)] the natural map of sheaves $\eta\colon f^* f_*(F)\to F$ is
surjective.
\end{description}
\end{prop}

\begin{pf}
Let $F$ be a perverse ideal sheaf on $Y$ and $E$ the 
corresponding perverse structure sheaf. Applying the
homology functor to the exact sequence (\ref{birthday})
and using Lemma \ref{perverse} shows that $F$ is a sheaf.
The functor $\R f_*$ is exact on the category
of perverse sheaves, so there is an exact sequence of sheaves
\[0\lra f_*(F)\lra \OO_X\lra \R f_*(E)\lra 0.\]
It follows that $f_*(F)$ is an ideal sheaf on $X$.

Let $A$, $B$ and $C$ denote the kernel, cokernel and image of the map
$\eta$ in the category of sheaves on $Y$.
Thus we have a pair of short exact sequences fitting
into a diagram
\begin{equation*}
 \setlength{\unitlength}{36pt}
 \begin{picture}(6.2,1.2)(0,0)
 \put(2.95,0){\object{$C$}}
 \put(0,0.75){\object{$0$}}
 \put(1,0.75){\object{$A$}}
 \put(2.3,0.75){\object{$f^* f_* F$}}
 \put(3.8,0.75){\object{$F$}}
 \put(4.95,0.75){\object{$B$}}
 \put(5.95,0.75){\object{$0.$}}
 \put(2.35,0.45){\vector(1,-1){0.3}}
 \put(3.3,0.15){\vector(1, 1){0.3}}
 \put(0.25,0.75){\vector(1,0){0.5}}
 \put(1.25,0.75){\vector(1,0){0.5}}
 \put(2.85,0.75){\vector(1,0){0.6}}
 \put(4.15,0.75){\vector(1,0){0.5}}
 \put(5.2,0.75){\vector(1,0){0.5}}
 \put(3.15,0.75){\elabel{$\eta$}}
 \end{picture}
\end{equation*}
The spectral sequence of Lemma \ref{easy} gives an exact sequence
\[0\lra \R^1 f_*(\L_1 f^* f_*(F))\lra f_*(F)\lra f_*(f^* f_*(F))\lra
0.\]
together with the fact that $\R^1 f_*(f^* f_*(F))=0$.
Since $f_*(F)$ is torsion-free this implies that
$f_*(f^* f_*(F))=f_*(F)$.

The morphism $f$ is birational, so $\eta$ is generically an isomorphism
and  $A$ and $B$ are torsion sheaves. Applying $f_*$ to the exact
sequences above shows that $f_*(F)$ injects into $f_*(C)$ and also
$f_*(C)$ injects into $f_*(F)$. It follows that $f_*(C)=f_*(F)$ and so
$f_*(B)=0$. Since $F$ is perverse, $\R^1 f_*(F)=0$, so $\R f_*(B)=0$
and hence by Lemma \ref{perverse}, $B=0$, that is, $\eta$ is surjective.

For the converse, suppose $F$ is a sheaf on $Y$ satisfying our two
conditions. There is an exact sequence of sheaves
\begin{equation}
\label{lednock}
0\lra A\lra f^* f_*(F)\lRa{\eta} F\lra 0.
\end{equation}
It follows that $\R^1 f_*(F)=0$ and $\Hom_Y(F,C)=0$ for any sheaf $C$
on $Y$ satisfying $\R f_*(C)=0$, hence, by Lemma \ref{perverse}, $F$
is a perverse sheaf on $Y$.
Since $f$ is birational, $\eta$
is generically an isomorphism, so $A$ is a torsion sheaf and
\[
\Hom_Y(F,\OO_Y)=\Hom_Y(f^* f_*(F),\OO_Y)=\Hom_X(f_*(F),\OO_X)
\]
which is non-zero because $f_*(F)$ is an ideal sheaf.
Take a non-zero morphism $F\to\OO_Y$ and form a triangle
\[ F\lra \OO_Y\lra E\lra F[1]. \]
It will be enough to show that $E\in\PC(Y/X)$. Applying the homology
functor gives a long exact sequence of sheaves
\[ 0\lra H_1(E)\lra F\lra \OO_Y\lra H_0(E)\lra 0.\]
It follows that  $\R^1 f_*( H_0(E))=0$ and $\Hom_Y(H_0(E),C)=0$ for
any sheaf $C$ on $Y$ satisfying $\R f_*(C)=0$. Furthermore $f_*(F)$ is
torsion-free so $\R^0 f_* (H_1(E))=0$.
Applying Lemma \ref{perverse} completes the proof.
\end{pf}

\begin{lemma}
Any perverse ideal sheaf $F$ on $Y$ is simple, that is
$\Hom_Y(F,F)=\C$.
\end{lemma}

\begin{pf}
The ideal sheaf $f_*(F)$ is simple, so $\Hom_Y(f^* f_*(F),F)=\C$.
Applying the functor $\Hom_Y(-,F)$ to the sequence (\ref{lednock})
gives the result.
\end{pf}

\medskip

We shall now construct a fine moduli space of perverse ideal sheaves.
To do this we mimic C. Simpson's proof \cite[Section 1]{Si}
of the existence of moduli spaces of semistable sheaves
on projective schemes.

Consider the
following special
case of Simpson's construction. Take a sufficiently ample
line bundle $\OO_Y(1)$ on $Y$ and a suitable vector space $V$.
Then the stable points for the action of
the group $\SL(V)$ acting on the Quot scheme parameterising quotients
of $V\tensor \OO_Y(-1)$ with
rank one and trivial determinant
are just the points corresponding to ideal sheaves on $Y$.

We shall show that if one takes a sufficiently ample line bundle
$\OO_X(1)$ on $X$ and replaces $\OO_Y(-1)$ by $f^* \OO_X(-1)$ in the
above construction, then the stable points are precisely the points
corresponding to perverse ideal sheaves on $Y$.

Fix a numerical equivalence class $(\gamma)$ and let $F$ denote a
perverse ideal sheaf on $Y$ in this class.

Rank one, torsion-free sheaves in a given numerical equivalence class
form a bounded family, so we may choose
$\OO_X(1)$ so that for any torsion-free sheaf $A$ on $X$ in the
same numerical equivalence class as $f_*(F)$, the sheaf $A\tensor\OO_X(1)$
is generated by its global sections and satisfies
$H^i(X,A\tensor\OO_X(1))=0$ for all $i>0$.

Put $\W=f^*(\OO_X(-1))$ and let $V$ be the vector space
\[V=\Hom_Y(\W,F)=\Hom_X(\OO_X(-1),f_*(F)).\]
Let $\Quot$ denote the Quot scheme parameterising
quotients of the vector bundle $V\tensor\W$ in the numerical
equivalence class $(\gamma)$.

\begin{lemma}
There is a locally-closed subscheme $\U\subset\Quot$ 
parameterising quotients $V\tensor\W\onto F$
for which $f_*(F)$ is
an ideal sheaf and for which the natural map
$\alpha_F\colon V\to \Hom_Y(\W,F)$
is an isomorphism.
\end{lemma}

\begin{pf}
By Proposition \ref{basechange} there is an open subscheme of $\Quot$
parameterising quotients for which $f_*(F)$ is torsion-free.
There is also a closed subscheme over which there is a non-zero map
$f_*(F)\to\OO_X$. These two conditions are equivalent to
$f_*(F)$ being an ideal sheaf. The condition on the map
$\alpha_F$ is clearly
open.
\end{pf}

The group $\SL(V)$ acts on $\Quot$ preserving the subscheme $\U$.
Define $\Quot^0$ to be the closure of $\U$ in $\Quot$.
Fix a very ample line bundle $\OO_Y(1)$ on $Y$.
For sufficiently large integers $m$, there is a closed embedding of
$\Quot^0$ in a Grassmannian, which sends the quotient
$V\tensor\W\onto F$ to the quotient of vector spaces
\[\Hom_Y(\OO_Y(-m),V\tensor\W)\onto \Hom_Y(\OO_Y(-m),F).\]
This embedding is $\SL(V)$-equivariant. By pulling back the
natural $\SL(V)$-equivariant polarisation of the Grassmannian
one obtains an $\SL(V)$-equivariant polarisation $\Li(m)$ of $\Quot^0$.

\begin{prop}
For all sufficiently large integers $m$ the following
result holds. The semistable points for the action of $\SL(V)$ on the
projective scheme $\Quot^0$ with respect to the polarisation $\Li(m)$
are precisely the points of the open subset $\U$.
Furthermore all these points are properly stable.
\end{prop}

\begin{pf}
Let $m$ be larger than the integer $M$ given by \cite[Lemma 1.15]{Si}
and take a point $p\in \U$ corresponding to a quotient
$V\tensor\W\onto F$. Let $H\subset V$ be a proper non-zero subspace
and let $G\subset F$ be the subsheaf generated by $H\tensor\W$.
Since $\alpha_F$ is an isomorphism, $G$ is non-zero.
But $G$ is a quotient of $H\tensor\W$, so $f_*(G)$ is non-zero,
and since $f_*(F)$ is torsion-free this implies that $G$ has rank one.
Thus the leading coefficient of the Hilbert polynomial $P_G$ is the
same as that of $P_F$.

Since the class of subsheaves of
$F$ which are quotients of $V\tensor\W$ is bounded,
the set of possible Hilbert polynomials of $G$ is finite, so
we may assume that the inequality
\[\frac{\dim(H)}{P_G(m)}<\frac{\dim(V)}{P_F(m)}\]
holds. Lemma \cite[1.15]{Si}
then implies that the point $p$ is stable with repect to $\Li(m)$.

For the converse, let $m$ be larger than the integer $M$
given by \cite[Lemma 1.16]{Si}
and take a point $p\in\Quot^0$ which is semistable
with repect to $\Li(m)$.
Let $V\tensor\W\onto F$ be the corresponding quotient.
By \cite[Lemma 1.16]{Si}, for any non-zero subspace $H\subset V$,
the subsheaf generated by $H\tensor\W$ has positive rank.
It follows that the map $\alpha_F$ is an isomorphism.

Let $T$ be the torsion subsheaf of $f_*(F)$ and put $Q=f_*(F)/T$.
By assumption $f_*(F)$ is a degeneration of torsion-free sheaves,
so by \cite[Lemma 1.17]{Si} there exists a torsion-free sheaf $A$
with the same numerical invariants as $f_*(F)$ and an inclusion
$Q\subset A$.

By our choice of $\OO_X(1)$ the vector space $\Hom_X(\OO_X(-1),A)$
has the same dimension as $V$ and the
sheaf $A\tensor\OO_X(1)$ is generated by its global
sections.
Since $Q$ has rank one, the remark after
\cite[Lemma 1.16]{Si} shows that the natural map $V\to\Hom_Y(\W,f^* Q)$
is injective.
It follows that $Q=A$. Then $T$ is numerically trivial,
so $T=0$ and $f_*(F)$ is torsion-free. By semi-continuity there is
a non-zero map $f_*(F)\to\OO_X$ so $f_*(F)$ is an ideal sheaf.
\end{pf}

\begin{thm}
\label{moduli}
The functor which assigns to a scheme $S$ the set of equivalence
classes of families over $S$
of perverse ideal sheaves in a given numerical equivalence class
$(\gamma)$ is representable by a projective scheme $\M_{PI}(Y/X;\gamma)$.
\end{thm}

\begin{pf}
Let $F$ be a perverse ideal sheaf on $Y$ in the numerical
equivalence class $(\gamma)$.
By assumption $f_*(F)\tensor\OO_X(1)$ is generated
by its global sections so there exists a surjection
$V\tensor\OO_X(-1)\onto f_*(F)$. Pulling back and using
Proposition \ref{id} shows that there is a
surjection $V\tensor\W\onto F$ and hence a point of $\U$ for 
which the corresponding quotient is $F$.

Conversely, the argument of Lemma \ref{perverse} shows that
a quotient $V\tensor\W\onto F$ corresponding to a point of $\U$
is a perverse ideal sheaf.
Exactly as in \cite[Theorem 1.21]{Si} we can conclude that there is a coarse
moduli space $\M_{PI}(Y/X;\gamma)$
for perverse ideal sheaves in the numerical
equivalence class $(\gamma)$, and that a universal sheaf exists
locally in the {\'e}tale topology on $\M$.

All perverse ideal sheaves $F$ have rank one, so
$\eu(F,\OO_y)=1$ for any point $y\in Y$. If $y$ is a non-singular
point then $\OO_y$ has a finite locally-free resolution, so
the integers $\eu(F\tensor L),$ as $L$ ranges over all locally-free sheaves on
$Y$, have no common factor. Since the sheaves $F$ are simple,
an argument of S. Mukai \cite[Theorem A.6]{Mu} shows that one can
patch the local universal sheaves to obtain a universal sheaf on
$\M_{PI}(Y/X;\gamma)$. This completes the proof.
\end{pf}


\section{Perverse Hilbert schemes}

In this section we complete the proof of Theorem \ref{rep}.
To do this we construct a perverse Hilbert scheme $\pHilb(Y/X)$
parameterising quotients of $\OO_Y$ in the category $\PC(Y/X)$.
As before $f\colon Y\to X$ denotes a birational morphism of
projective varieties of relative dimension one and satisfying
$\R f_*(\OO_Y)=\OO_X$.

\medskip

If $S$ is a scheme, the $S$-valued points of the
Hilbert scheme of a variety $Y$ consist of triangles
$\F\lra\OO_{S\times Y}\lra\E\lra \F[1]$
in $\D(S\times Y)$ such that $\E$ and $\F$ are families of sheaves on $Y$
over $S$. Analagously we make the

\begin{defn}
Let $\pHilb(Y/X)$ be the functor which assigns to a scheme
$S$ the set of isomorphism classes of triangles
$\F\lra\OO_{S\times Y}\lra\E\lra \F[1]$ in $\D(S\times Y)$
such that $\F$ and $\E$ define families of perverse sheaves on $Y$
over $S$.
\end{defn}

Given a scheme $S$ we write
$f_S=\id_S\times f\colon S\times Y\to S\times X$.

\begin{lemma}
\label{bonza}
If $\F$ is a family  of perverse ideal sheaves on $Y$
over $S$, the map 
\[\fs\colon\Hom_{S\times Y}(\F,\OO_{S\times Y})\lra
\Hom_{S\times X}(\fs(\F),\OO_{S\times X})\]
is an isomorphism.
\end{lemma}

\begin{pf}
Since perverse ideal sheaves are sheaves, the object $\F$ is a sheaf
on $S\times Y$, flat over $S$. Consider
the natural map of sheaves
$\eta\colon f_S^* \fs (\F)\to \F$.
For each point $s\in S$, the map $\eta_s=\L j_s^*(\eta)$
is just the natural map $f^* f_* (\F_s)\to\F_s$ which is surjective
by Proposition \ref{id}. It follows that $\eta$ is surjective.

Let $K$ be the kernel of $\eta$.
It will be enough
to show that there are no non-zero maps
$K\to\OO_{S\times Y}$. For this we may assume that $S$ is
affine.
Let $p\colon S\times Y\to Y$ be the projection.
Since $f$ is birational, $p_*(K)$ is a torsion sheaf. But
$p_*(\OO_{S\times Y})$ is torsion-free, so any
map $p_*(K)\to p_*(\OO_{S\times Y})$ is zero.
\end{pf}  

The functor  $\pHilb(Y/X)$ decomposes into components which parameterise
ideal sheaves in a given numerical
equivalence class $(\gamma)$:
\[\pHilb(Y/X)=\coprod_{(\gamma)}\pHilb(Y/X;\gamma)\]

\begin{thm}
\label{hilb}
For each numerical equivalence class $(\gamma)$ the corresponding functor
$\pHilb(Y/X;\gamma)$ is representable by a projective scheme.
\end{thm}

\begin{pf}
By Theorem \ref{moduli} there is a
fine moduli space $\M_{PI}(Y/X)$ for perverse ideal sheaves on $Y$.
An $S$-valued point of $\M_{PI}(Y/X)$ is a family of perverse
ideal sheaves on $Y$ over $S$. By Proposition \ref{basechange},
applying the functor $\R f_*$ gives a family
of ideal sheaves on $X$ over $S$. This gives a morphism
\[\M_{PI}(Y/X)\lra\M_I(X),\]
where the scheme on the right is the moduli space of ideal sheaves
on $X$.

On the other hand, an $S$-valued point of the Hilbert
scheme  $\Hilb(X)$ determines a family of ideal sheaves on $X$ over $S$,
so there is a morphism
\[\Hilb(X)\lra\M_I(X),\]
which induces a bijection on closed points.
I claim that $\pHilb(Y)$ is represented by the fibre product
\begin{equation*}
\setlength{\unitlength}{54pt}
\begin{picture}(2.2,1.2)(0,0)
\put(0,0){\object{$\Hilb(X)$}}
\put(2,0){\object{$\M_I(X)$}}
\put(0,1){\object{$\pHilb(Y/X)$}}
\put(2,1){\object{$\M_{PI}(Y/X)$}}
\put(0.75,.45){$\qed$}
\put(0.8,0){\vector(1,0){.5}}
\put(0,0.75){\vector(0,-1){.5}}
\put(0.8,1){\vector(1,0){.5}}
\put(2,0.75){\vector(0,-1){.5}}
\end{picture}
\end{equation*}

Indeed, an $S$-valued point of the functor $\pHilb(Y/X)$ is a triangle
of objects
\[\F\lra\OO_{S\times Y}\lra\E\lra \F[1]\]
of $\D(S\times Y)$,
each of which is
a family of perverse sheaves
on $Y$ over $S$. Applying the functor $\R \fs$ gives a short exact sequence
of sheaves on $S\times X$, each of which is a family of
sheaves on $X$ by Proposition \ref{basechange}. This defines an
$S$-valued point of $\Hilb(X)$. The family $\F$ defines an
$S$-valued point of $\M_{PI}(Y/X)$ so we get an $S$-valued point
of the fibre product.

Conversely, an $S$-valued point of the fibre product
gives a family of perverse ideal sheaves $\F$ on $Y$ over $S$ together
with a short exact sequence
\[ 0\lra \fs(\F)\lra \OO_{S\times X}\lra \G\lra 0\]
of $S$-flat sheaves on $S\times X$.
By Lemma \ref{bonza} this uniquely determines a triangle
\[\F\lRa{\alpha} \OO_{S\times Y}\lra \E\lra \F[1]\]
in $\D(S\times Y)$ with $\G=\R \fs(\E)$. This is an $S$-valued point of
the functor $\pHilb(Y/X)$ providing $\E$ defines a family of perverse
sheaves on $Y$ over $S$.

Applying the functor $\L j_s^*$ and noting
that, as in Lemma \ref{id},
any non-zero morphism $\F_s\to\OO_Y$ is an injection in
$\PC(Y/X)$, it is enough to show that $\L j_s^*(\alpha)$
is non-zero for all $s\in S$.
But if $\L j_s^*(\alpha)$ vanishes then
\[ H^{-1}(\L j_s^*(\E))=H^0(\F_s) \]
so $\R f_* (\L j_s^*(\E))$ cannot
be a sheaf. But $\G$ is flat over $S$, so this
contradicts the argument of Proposition \ref{basechange}.
\end{pf}

\medskip

\noindent {\bf Proof of Theorem \ref{rep}.}
Let $(\gamma)$ denote the numerical equivalence class of the ideal
sheaf $\I_y$ of a point $y\in
Y$. I claim that the scheme $\pHilb(Y/X;\gamma)$ is a fine moduli space
for perverse point sheaves on $Y$.

An $S$-valued point of $\pHilb(Y/X;\gamma)$ certainly determines
a family of perverse point sheaves on $Y$ over $S$. For the converse
suppose $\E$ is a family of perverse point sheaves on $Y$ over $S$.
The object $\G=\R \fs(\E)$ is a family of points on $X$ over
$S$ and hence, up to a twist by a line bundle from $S$, is the
structure sheaf of the graph of a morphism $S\to X$.
Twisting $\E$ by the pullback of a line bundle from $S$
we may assume that there is a surjection
$\delta\colon\OO_{S\times X}\onto \G$
whose kernel is flat over $S$.

By adjunction, there is
a morphism $\beta\colon\OO_{S\times Y}\to \E$ such
that $\R \fs(\beta)=\delta$. Forming a triangle
\[\F\lra\OO_{S\times Y}\lRa{\beta}\E\lra \F[1]\]
gives an $S$-valued point of $\pHilb(Y/X;\gamma)$
providing $\F$ is a family of perverse ideal sheaves.
Applying the functor $\L j_s^*$ it will be enough to check that
$\L j_s^*(\beta)$ is nonzero for all $s\in S$.
But if $\L j_s^*(\beta)=0$ then
\[H^1(\L j_s^*(\F))=H^0(\E_s)\]
so $\R f_*(\L j_s^*(\F))$ cannot be a sheaf,
and this contradicts the argument of
Proposition \ref{basechange}, since $\R \fs(\F)$ is flat over $S$. This
completes the proof.
\qed


\bigskip

\noindent Department of Mathematics and Statistics,
The University of Edinburgh,
King's Buildings, Mayfield Road, Edinburgh, EH9 3JZ, UK.

\smallskip

\noindent email: {\tt tab@maths.ed.ac.uk}

\end{document}